\input amstex                                                           \documentstyle{amsppt}
\voffset=-3pc
\def\bC{\Bbb C}

\def\bR{\Bbb R}

\def\cQ{\Cal Q}
\def\cK{\Cal K}

\def\cS{\Cal S}

\def\cW{\Cal W}
\def\ep{\epsilon}

\def\ep{\epsilon}

\def\b1{\bold 1}

\def\tf{\widetilde f}

\def\im{\text{Im}}
\def\re{\text{Re}}
\def\m{\Bbb M_{n\text {sa}}}
\magnification=\magstep1
\parskip=6pt
\NoBlackBoxes
\topmatter
\title A Treatment of Strongly Operator convex Functions That Does Not Require
Any Knowledge of Operator Algebras
\endtitle  
\author Lawrence G.~Brown
\endauthor
\rightheadtext{Strongly Operator Convex Functions}
%\endtopmatter

\abstract{In [B1, Theorem 2.36] we proved the equivalence of six conditions on a continuous function $f$ on an interval. These conditions determine a subset of the set of operator convex functions, whose elements are called strongly operator convex. Two of the six conditions involve operator--algebraic semicontinuity theory, as given by C. Akemann and G. Pedersen in [AP], and the other four conditions do not involve operator algebras at all. Two of these conditions are operator inequalities, one is a global condition on $f$, and the fourth is an integral representation of $f$, stronger than the usual integral representation for operator convex functions. The purpose of this paper is to make the equivalence of these four conditions accessible to people who do not know operator algebra theory as well as to operator algebraists who do not know the semicontinuity theory. We also provide a similar treatment of one theorem from [B1] concerning (usual) operator convex functions. And in two final sections we give a somewhat tentative treatment of some other operator inequalities for strongly operator convex functions, and we give a differential criterion for strong operator convexity.}
\endabstract
\endtopmatter

\subhead 1. \ Introduction\endsubhead

A continuous real--valued function $f$ on an interval $I$ is called \underbar{operator monotone} if $h_1\le h_2$ and $\sigma(h_i)\subset I$ imply $f(h_1)\le f(h_2)$ and \underbar{operator convex} if $\sigma(h_1), \sigma(h_2)\subset I$ and $t\in [0,1]$ imply $f(th_1+(1-t)h_2)\le tf(h_1)+(1-t)f(h_2)$. Here $h_1$ and $h_2$ are in $B(H)_{\text{sa}}$, the set of self--adjoint bounded operators on a Hilbert space $H$, $\sigma(h_i)$ denotes the spectrum, and $f(h_i)$ is defined by the continuous functional calculus. Because of the assumed continuity of $f$ it is sufficient to verify either condition for finite dimensional $H$, and operator monotonicity or convexity on the interior of $I$ implies the same on all of $I$. The reader is referred to [D2] and [H] for further information on these topics and their history. C. Davis showed in [D1] that $f$ is operator convex if and only if

\noindent
(1)\   $pf(php)p\le p f(h)p$ for $h\in B(H)_{\text{sa}}$ with $\sigma(h)\subset I$ and $p$ a projection.\newline

And F. Hansen and G. Pedersen showed in [HP] that if $0\in I$ and $f(0)\le 0$, then $f$ is operator convex if and only if

\noindent
(2)\  $f(a^*ha)\le a^*f(h)a$ for $h\in B(H)_{\text{sa}}$ with $\sigma(h)\subset I$ and $||a||\le 1$.

\noindent
For comparative purposes we state [B1, Theorem 2.36], slightly rephrased, as well as the version that will be proved in this paper.

\proclaim{Theorem 0}(cf [B1, Theorem 2.36]) If $f$ is a continuous real--valued function on an interval $I$ containing $0$, then the following are equivalent.
\itemitem{(i)} If $h$ is a self--adjoint quasimultiplier of $E$ and $\sigma(H)\subset I$, then $f(h)$ is strongly lower semicontinuous.
\itemitem{(ii)} If $p, h\in B(H)_{\text{sa}}$ such that $p$ is a projection and $\sigma(h)\subset I$, then $pf(php)p\le f(h)$. 
\itemitem{(iii)} If $p, h\in B(H)_{\text{sa}}$ such that $0\le p\le \b1$ and $\sigma(h)\subset I$, then $f(php)\le f(h)+f(0)(\b1-p)$
\itemitem{(iv)} The condition in (i) holds if $E$ is replaced by an arbitrary $C^*$--algebra $A$.
\itemitem{(v)} Either $f=0$ or $f(x)>0$, $\forall x\in I$, and $-1/f$ is operator convex.
\itemitem{(vi)} $f$ has a representation

\noindent
(3) $ f(x)=c+\int_{r<I}{1\over x-r} d\mu_-(r)+\int_{r>I}{1\over r-x}d \mu_+(r)$, where $\mu_{\pm}$ are positive measures such that $\int{1\over 1+|r|} d\mu_{\pm}(r)<\infty$ and $c\ge 0$.
\endproclaim

\proclaim{Theorem 3.2} If $f$ is a continuous real--valued function on an interval $I$ containing 0, then the following are equivalent.
\itemitem{(i)} If $h\in \cQ$ and $\sigma (h)\subset I$, then $f(h)\in \cS$.
\itemitem{(ii)}Same as (ii) of Theorem 0.
\itemitem{(iii)}Same as (iii) of Theorem 0.
\itemitem{(iv)}Same as (v) of Theorem 0.
\itemitem{(v)}Same as (vi) of Theorem 0.
\endproclaim

In the above $E$ is a specific $C^*$--algebra whose definition will be given below (though we don't really need it), $\b1$ is the identity of $B(H)$, and $\cQ$
and $\cS$ will be defined below in concrete ways. Since Theorem 3.2 (i) and Theorem 0 (i) are in fact equivalent, we will provide references along with some of our preliminary results to the corresponding semicontinuity results, for the benefit of those readers who have some interest in the semicontinuity theory. However, the concepts and proofs below will use only some basic operator theory.

The requirement that 0 be in $I$  was discussed in [B1, Remark 2.37 (a)] and will similarly be discussed in Remark 3.3 (i) below. Both operator convexity and strong operator convexity are invariant under translation of the independent variable.

Another characterization of strong operator convexity is given in [B2, Theorem 4.8]. Since this seems to depend essentially on operator algebraic semicontinuity theory, it will not be further mentioned in this paper.

\subhead 2. \ Preliminaries\endsubhead

We will work with operators on the Hilbert space $\ell^2$. The set of compact operators on $\ell^2$ is denoted  by $\cK$, and $v\times w$, for $v,w$ vectors in $\ell^2$, is the operator $u\to (u,w)v$. If $h\in B(\ell^2)_{\text{sa}}$, we  can write uniquely $h=h_+-h_-$ where $h_+, h_-\ge 0$ and $h_+h_-=0$. Let $E$ be the set of norm convergent sequences in $\cK$, and denote by $E^{**}$ the set of bounded indexed collections $\{t_n\}_{1\le n\le\infty}$ with each $t_n$ in $B(\ell^2)$. Algebraic operations on $E$ and $E^{**}$, including the operation $t\mapsto t^*$, are defined componentwise; and for $k=(k_n)$ in $E$ or $t=\{t_n\}$ in $E^{**}$, $\| k\|=\sup\{\|k_n\|: \ 1\le n<\infty\}$ and $\|t\|=\sup\{\|t_n\|: \ 1\le n\le\infty\}$. It is not necessary to know that $E^{**}$ can be identified with the Banach space bidual of $E$, but it may be helpful to keep in mind that $E^{**}$ is a Banach algebra. We will denote by $\b1=\b1_{E^{**}}$ the elements $\{t_n\}$ of $E^{**}$ such that $t_n=\b1_{\ell^2}$ for $1\le n\le\infty$, and for $h=\{h_n\}$ in $E^{**}_{\text{sa}}$, $\sigma(h)$ denotes the set $(\cup_{1\le n\le\infty} \ \sigma(h_n))^-$ where $E^{**}_{\text{sa}}$ is the set of self--adjoint elements of $E^{**}$. If $f$ is a continuous function whose domain includes $\sigma(h)$, for $h=\{h_n\}$ in $E^{**}_{\text{sa}}$, then $f(h)$ denotes $\{f(h_n)\}_{1\le n\le\infty}$. Finally, if $h'=\{h'_n\}$ and $h^{''}=\{h^{''}
_n\}$ are two elements of $E^{**}_{\text{sa}}$, $h'\le h^{''}$ means $h'_n\le h^{''}_n$ for $1\le n\le\infty$. 

\proclaim{Lemma 2.1} If $h\in B(\ell^2)_{\text{sa}}$, then the following are equivalent.
\itemitem{(i)} $h_-\in\cK$.
\itemitem{(ii)} There is an increasing sequence $(k_n)$ in $\cK$ such that $k_n\to h$ weakly (equivalently, strongly). 
\itemitem{(iii)} There is $k\in\cK$ such that $h\ge k$.
\endproclaim

\demo{Proof}\ (i)$\Rightarrow$ (ii):\  Let $(p_n)$ be an increasing sequence of finite rank projections such that $p_n\to 1$ weakly and take $k_n=h_+^{1\over 2}p_nh_+^{1\over 2}-h_-$.

(ii)$\Rightarrow$ (iii):\  Let $k=k_1$.

(iii)$\Rightarrow$ (i):\  Since $h_+-h_-\ge k$, then $h_-\le h_+-k$. Therefore $h_-^3=h_-h_-h_-\le h_-h_+h_--h_-kh_-=-h_-kh_-$. The facts that $0\le h_-^3\le -h_-kh_-$ and $-h_-kh_-\in \cK$ imply $h_-^3\in \cK$, whence $h_-\in\cK$.
\enddemo

\proclaim{Lemma 2.2} (cf. [B1, 5.3]) Let $h$ in $B(\ell^2)_{\text{sa}}$ be the weak limit of an increasing sequence $(k_n)$ with each $k_n$ in $\cK$, let $k\le h$ for $k\in\cK$, and let $\ep>0$. Then $k\le k_n+\ep\b1$ for $n$ sufficiently large.
\endproclaim

\demo{Proof}\ If this is false, then there are unit vectors $v_n$ such that $(k v_n, v_n)>(k_nv_n, v_n)+\ep$, $\forall n$. Choose a subsequence $(v_{n_i})$ which converges weakly to a vector $v$. Since $(k_nv, v)\to (hv,v)\ge (kv, v)$, then $(k_n v, v)>(kv,v)-{\ep\over 3}$ for $n$ sufficiently large.  Choose one such $n$. Since $kv_{n_i}\to kv$ in norm and $k_n v_{n_i}\to k_nv$ in norm, then 
$(kv_{n_i}, v_{n_i})\to (kv,v)$ and $(k_n v_{n_i}, v_{n_i})\to (k_n v, v)$.  Therefore for $i$ sufficiently large, $n_i\ge n$, $|(kv_{n_i}, v_{n_i})-(kv, v)|< {\ep\over 3}$, and $|(k_nv_{n_i}, v_{n_i})-(k_nv, v)|< {\ep\over 3}$. Choose one such $i$. Thus $(kv, v)>(kv_{n_i}, v_{n_i})-{\ep\over 3}> (k_{n_i} v_{n_i}, v_{n_i})+{2\ep\over 3}\ge (k_n v_{n_i}, v_{n_i})+{2\ep\over 3}>(k_n v, v)+{\ep\over 3}$,  a contradiction. 
\enddemo

\proclaim{Definition 2.3} We denote by $\cS$ the set of elements $h=\{h_n\}$ in $E^{**}_{\text{sa}}$ such that:
\itemitem{(i)} $h_n$ satisfies the conditions in Lemma 2.1 for $1\le n\le\infty$, and
\itemitem{(ii)} If $k\in \cK$, $k\le h_\infty$, and $\ep>0$, then $k\le h_n+\ep\b1$ for $n$ sufficiently large. 
\endproclaim

Note that for $\lambda\in\bR$, $\lambda\b1\in \cS$ if and only if $\lambda\ge 0$.

\proclaim{Corollary 2.4} (cf. [B1, Remark (i) after 5.13]) Let $h=\{h_n\}$ be an element of $E^{**}_{\text{sa}}$ which satisfies 2.3 (i), and let $(k_m)$ be an increasing sequence in $\cK$ which converges weakly to $h_\infty$. If for each $m$ and each $\ep>0$, we have $k_m\le h_n+\ep\b1$ for $n$ sufficiently large, then $h\in\cS$. 
\endproclaim

\demo{Proof}\ Given $k$ in $\cK$ with $k\le h_\infty$ and $\ep>0$, apply Lemma 2.2 with $\ep/2$ in place of $\ep$.
\enddemo

\proclaim{Lemma 2.5}\ If $h=\{h_n\}$ is in $E^{**}_{\text{sa}}$, then the following are equivalent.
\itemitem{(i)}For each vector $v$, $(h_\infty v,v)\le \lim\inf(h_nv,v)$.
\itemitem{(ii)}For each weak cluster point $h'$ of the sequence $(h_n)$, $h_\infty\le h'$.
\itemitem{(iii)} For each finite rank projection $p$ and each $\ep>0$, $ph_\infty p\le ph_np+\ep p$ for $n$ sufficiently large.
\endproclaim

\demo{Proof}\  (i) $\Rightarrow$ (ii): \ For each vector $v$, $(h'v,v)$ is a cluster point of $((h_n v,v))$. Therefore $(h'v,v)\ge\lim\inf(h_n v,v)\ge (h_\infty v,v)$, whence $h'\ge h_\infty$.

(ii) $\Rightarrow$ (iii): \ If false, there is a subsequence $(h_{n_i})$ such that the relation $ph_\infty p\le ph_{n_i} p +\ep p$ is false, $\forall i$. Passing to a further subsequence, we may assume $h_{n_i}\to h'$ weakly for some $h'$. Then $ph_{n_i}p\to ph'p$ in norm. Therefore for $i$ sufficiently large, $ph_{n_i}p\ge ph'p-\ep p\ge ph_\infty p-\ep p$, a contradiction.

(iii) $\Rightarrow$ (i):\ For a unit vector $v$, let $p$ be the rank one projection $v\times v$. Since $ph_\infty p=(h_\infty v,v)p$ and $ph_np=(h_n v, v)p$, the given relation implies that $\forall \ep>0$, we have $(h_\infty v,v)\le (h_n v,v)+\ep$ for $n$ sufficiently large. Therefore $(h_\infty v,v)\le \lim\inf(h_n v,v)$.
\enddemo

\proclaim{Definition 2.6}\ Denote by $\cW$ the set of $h$ in $E^{**}_{\text{sa}}$ satisfying the conditions in Lemma 2.5, and denote by $\cQ$ the set of $h=\{h_n\}$ in $E^{**}_{\text{sa}}$ such that $h_n\to h_\infty$ weakly. Thus $h\in\cQ$ if and only if $h\in\cW$ and $-h\in\cW$. Note that $\lambda\b1\in\cQ\subset \cW, \forall \lambda\in\bR$. 
\endproclaim

If $t$ is in $B(\ell^2)$ and $h=\{h_n\}$ is in $E^{**}_{\text{sa}}$, then $t^*ht$ denotes the element $\{t^* h_n t\}$ of $E^{**}_{\text{sa}}$

\proclaim{Proposition 2.7} The sets $\cS$ and $\cW$ are closed in the norm topology and are closed under addition, multiplication by non-negative scalers, and the operation $h\mapsto t^* ht$, $t\in B(\ell^2)$. Also $\cS\subset \cW$.
\endproclaim

\demo{Proof}\ It follows easily from 2.5 (i) or 2.5 (iii) that $\cW$ is norm closed. Suppose $h^{(m)}\in\cS$ for $m=1,2,\dots$, and $h^{(m)}\to h$ in the norm of $E^{**}$. Since the map $t\mapsto t_-$ is norm continuous on $B(\ell^2)_{\text{sa}}$, it is clear that $h_n$ satisfies 2.1 (i) for $1\le n\le\infty$. Now let $p$ be a finite rank projection, $k=(h_{\infty +})^{1\over 2} p(h_{\infty +})^{1\over 2}- h_{\infty -}$ and $k^{(m)}=(h^{(m)}_{\infty +})^{1\over 2} p (h^{(m)}_{\infty +})^{1\over 2} -h^{(m)}_{\infty -}$. Then $k, k^{(m)}\in \cK$, $k\le h_\infty$, $k^{(m)}\le h^{(m)}_\infty$, and  $k^{(m)}\to k$ in norm. Let $\ep>0$ and choose an $m$ such that $\| h^{(m)}-h\|<{\ep\over 3}$ and $\| k^{(m)}-k\|<{\ep\over 3}$. Then $\exists N$ such that $n>N\Rightarrow k^{(m)}\le h^{(m)}_n+{\ep\over 3}\b1$. Then for $n>N$, $k\le k^{(m)}+{\ep\over 3}\b1\le h^{(m)}_n+{2\ep\over 3}\b1\le h_n+\ep\b1$. Thus by Corollary 2.4 and the proof of (i)$\Rightarrow$ (ii) in Lemma 2.1, we conclude that $h\in \cS$.

It is obvious that $\cS$ and $\cW$ are closed under multiplication by non-negative scalars. Let $h', h''\in\cS$ and $h=h'+h''$. Clearly $h_n$ satisfies 2.1 (iii) for $1\le n\le\infty$. If $(k'_m)$ and $(k''_m)$ are increasing sequences in $\cK$ such that $k'_m\to h'_\infty$ and $k''_m\to h''_\infty$ weakly, then apply corollary 2.4 with $k_m=k'_m+k''_m$ to conclude that $h\in\cS$. The situation is similar for the operator $h\mapsto t^* ht$. It is obvious for $\cW$ (note that ($t^* h_n t v, v)=(h_n tv, tv)$), and for $\cS$ we use the sequence $(t^* k_m t)$, which increases to $t^*h_\infty t$ if $(k_m)$ increases to $h_\infty$. If $k_m\le h_n+\ep\b1$, then $t^* k_m t\le t^*h_n t+\ep\|t\|^2\b1$. 

Finally let $h\in\cS$ and choose a vector $v$. If $k\in\cK$ and $k\le h_\infty$, then the fact that $\forall \ep>0$, $k\le h_n+\ep\b1$ for $n$ sufficiently large implies that $(kv, v)\le \lim\inf(h_n v, v)+\ep\|v\|^2$. Since $\ep$ is arbitrary, this implies $(kv, v)\le\lim\inf (h_n v, v)$. And since $k$ is arbitrary and $h_\infty$ satisfies 2.1 (ii), this implies $(h_\infty v, v)\le \lim\inf (h_n v, v)$. 
\enddemo

\proclaim{Proposition 2.8} (cf. [AP Proposition 3.5], which is slightly rephrased in [B1, Proposition 2.1 (a)]). Assume $h\in E^{**}_{\text{sa}}$ and $h\ge \eta \b1$ for some $\eta >0$. Then $h\in \cS$ if and only if $-h^{-1}\in \cW$. 
\endproclaim

\demo{Proof}\ By replacing $h$ with $h_\infty^{-{1\over 2}}h h_\infty^{-{1\over 2}}$, we reduce to the case $h_\infty=\b1$.

Now assume $h\in\cS$, $p$ is a finite rank projection, and $\ep>0$. Choose $\delta >0$ such that $p(h_n+2\delta\b1)^{-1}p\ge p h^{-1}_n p-\ep p$, $\forall n$. Then $\exists N$ such that $p\le h_n+\delta\b1$ for $n> N$. Therefore $p+\delta\b1\le h_n+2\delta\b1$ for $n>N$, whence $(1+\delta)^{-1}p+\delta^{-1}(\b1 -p)=(p+\delta\b1)^{-1}\ge (h_n+2\delta\b1)^{-1}$ for $n>N$. Therefore $p\ge (1+\delta)^{-1}p\ge p(h_n+2\delta\b1)^{-1} p\ge ph^{-1}_n p-\ep p$ for $n> N$. Thus $-h^{-1}$ satisfies 2.5 (iii).

Next assume $-h^{-1}\in\cW$, $p$ is a finite rank projection, and $\ep>0$. Choose $\delta>0$ such that $(1+2\delta)h_n\le h_n+\ep\b1$, $\forall n$. 
Then $\exists N$ such that $ph^{-1}_n p\le p+\delta p$ for $n>N$. It follows that for some $\lambda >0$, $h^{-1}_n\le p+2\delta p+\lambda (\b1-p)$ for $n>N$. To see this, it is convenient to represent elements of $B(\ell^2)_{\text sa}$ by $2\times 2$ matrices $\pmatrix a & b\\ b^* & c\endpmatrix$, where $a\in p B(\ell^2)p$, $b\in pB(\ell^2)(\b1-p)$, etc. If $a\ge \eta_1p$ and $c\ge \eta_2(\b1-p)$ for $\eta_1, \eta_2>0$, then this matrix is positive if and only if $\|a^{-\frac 12} b c^{\frac 12}\|\le 1$. Now we have $h_n\ge (p+2\delta p+\lambda(\b1-p))^{-1}=(1+2\delta)^{-1}p+\lambda^{-1}(\b1-p)\ge (1+2\delta)^{-1}p$ for $n>N$. Thus $p\le h_n+\ep\b1$ for $n>N$. Since there is a sequence of finite rank projections which increases to $\b1$ and since $h\ge 0$, it follows that $h\in \cS$
\enddemo

\proclaim{Lemma 2.9}\ If $h\in\cQ$, then $h^2\in\cW$
\endproclaim

\demo{Proof}\ If $v\in\ell^2$, then $h_nv\to h_\infty v$ weakly. Also $(h^2_nv,v)=\|h_n v\|^2$ and $(h_\infty^2 v, v)=\|h_\infty v\|^2$. It is well known that, for vectors in a Hilbert space, $w_n\to w$ weakly implies $\|w\|\le\lim\inf \|w_n\|$. 
\enddemo

\subhead 3. \ Main results\endsubhead

One direction of the equivalence in the next theorem is needed for the proof of Theorem 3.2. We prove both directions because of the intrinsic interest.

\proclaim{Theorem 3.1}(cf [B1, Propositions 2.34 and 2.35(b)]) (a) If $f$ is a continuous real--valued function on a compact interval $[a,b]$, then $f$ is operator convex if and only if whenever $h_n\to h$ weakly where $h_n\in B(\ell^2)_{\text sa}$ and $\sigma(h_n)\subset [a,b]$, $\forall n$, and $v\in\ell^2$, then $(f(h)v, v)\le\lim\inf (f(h_n)v, v)$.

(b) Equivalently, if $f$ is a continuous real--valued function on an interval $I$, then $f$ is operator convex if and only if $h\in\cQ$ and $\sigma(h)\subset I$ imply $f(h)\in\cW$.
\endproclaim

\example{Remarks} (i) As is well known, the strong convergence of a sequence $(h_n)$ to $h$ in $B(\ell^2)_{\text sa}$ implies that $f(h_n)\to f(h)$ strongly for any continuous function $f$, but there is no similar implication for weak convergence. Version (a) says that operator convexity is characterized by the fact that $h_n\to h$ weakly implies ``half''  of what is needed to conclude that $f(h_n)\to f(h)$ weakly. In fact, if the operator convex function $f$ is non-linear, it is impossible that $f(h_n)\to f(h)$ weakly unless $h_n\to h$ strongly. This follows from [B1, Proposition 2.59 (a)]. The original plan for this paper was to include a non-operator algebraic proof of this, but it turns out that the result has nothing to do with operator convexity. If $f$ is merely a continuous strictly convex function, $h_n\to h$ strongly. This is proved in [B3] in an elementary way. (Of course, every non-linear operator convex function is strictly convex.)

(ii) The forward implication in version (a) can be strengthened by replacing $\ell^2$ with an arbitrary (possibly non-separable) Hilbert space and replacing the sequence $(h_n)$ with a (necessarily bounded) net. Essentially the same proof works, or the stronger version can be deduced from the version stated.
\endexample

\demo{Proof of Theorem 3.1}\  We prove version (b). We reduce to the case $0\in I$ and $f(0)=0$ by replacing $f$ with $f(\cdot+x_0)-f(x_0)$ for some $x_0$ in $I$. This does not affect either half of the claimed equivalence.

If $f$ is operator convex, then $f$ has a representation
$$
f(x)=ax^2+bx+c+\int_{r<I}{(x-x_0)^2\over (x-r)(x_0-r)^2} d_{\mu_-}(r)+\int_{r>I}{(x-x_0)^2\over (r-x)(r-x_0)^2} d_{\mu_+}(r),\leqno(4)
$$
where $x_0$ can be any interior point of $I$, $a\ge 0$, $b, c\in\bR$, and $\mu_{\pm}$ are positive measures such that $\int{1\over 1+|r|^3} d_{\mu_{\pm}}(r)<\infty$. If $I$ contains one or both of its endpoints, then convergence of (4) at such endpoint(s) imposes an additional condition on $\mu_\pm$. If $h\in\cQ$ and $\sigma(h)\subset I$, then $f(h)$ is obtained by substituting $h$ for $x$ in (4), thus obtaining a Bochner integral. (Note that the integrands in (4) give a continuous function from $\bR\setminus I$ to the Banach space $E^{**}$.) Because of the properties of $\cW$ proved in Proposition 2.7, it is enough to show that each value of the integrand and each term $ah^2,  bh, c\b1$ is in $\cW$. Now Proposition 2.8 implies that $(r\b1-h)^{-1}$, for $r>I$, and $(h-r\b1)^{-1}$, for $r<I$, are in $\cS\subset \cW$. Also the integrands in (4) are obtained from $1/(r-x)$ or $1/(x-r)$ by subtracting its first degree Taylor polynomial at $x=x_0$. Since the linear terms are in $\cQ\subset \cW$, and since Lemma 2.9 covers the $ah^2$ term, we conclude that $f(h)\in \cW$.

Now assume that $h\in \cQ$ and $\sigma(h)\subset I$ imply $f(h)\in\cW$. We will prove that $f$ is operator convex by proving (1), and we begin with a matrix version. For natural numbers $k,l$ consider $(k+l)\times(k+l)$ self-adjoint matrices
$$t=\pmatrix a & b\\ b^* & c\endpmatrix\quad\text{and}\quad p=\pmatrix \b1_k &0\\0 & 0\endpmatrix,$$
where $a$ is $k\times k$, $b$ is $k\times l$, etc., and $\sigma(t)\subset I$. Let $f(t)=\pmatrix a' & b'\\ b^{'*} & c'\endpmatrix$. Then the desired relation, $pf(ptp)p\le pf(t)p$, amounts to $f(a)\le a'$. Let $e_1, e_2,\dots$ be the standard orthonormal basis vectors for $\ell^2$, and define $h=\{h_n\}$ by $h_n=\sum a_{ij}e_i\times e_j+\sum b_{ii'}e_i\times e_{n+k+i'}+\sum \Bar{b}_{ii'} e_{n+k+i'}\times e_i+\sum c_{i'j'}e_{n+k+i'} \times e_{n+k+j'}$, for $n<\infty$, and $h_\infty=\sum a_{ij}e_i\times e_j$. Here $i,j=1, \dots, k$ and $i', j'=1, \dots, l$. Then $\sigma(h)=\sigma(t)\cup\sigma(a)\cup\{0\}\subset I$, and $h\in\cQ$. So $f(h)\in\cW$. If $f(h)=\{s_n\}_{1\le n\le\infty}$, then for finite $n$, $s_n$ has a similar formula to $h_n$ with $a, b, c$ replaced by $a', b', c'$. Thus $(s_n)$ converges weakly to $\sum a'_{ij}e_i \times e_j$, and our desired relation follows from $s_\infty\le \lim s_n$.

The general case, where $t\in B(H)_{\text {sa}}$ and $p$ is a projection in $B(H)$, follows by a standard argument: Let $(p_i)$ and $(q_i)$ be nets of finite rank projections such  that $p_i\le p$, $q_i\le \b1-p$, $p_i\to p$, and $q_i\to \b1-p$, with convergence in the strong operator topology. Then $(p_i+q_i)t(p_i+q_i)\to t$ strongly and $\sigma((p_i+q_i)t(p_i+q_i))\subset I$. So it is enough to prove (1) for $p_i$ and $(p_i+q_i)t(p_i+q_i)$, and this follows from the matrix version.
\enddemo

\proclaim{Theorem 3.2}\ (cf. [B1, Theorem 2.36]) If $f$ is a continuous real--valued function on an interval $I$ containing 0, then the following are equivalent.
\endproclaim

\itemitem{(i)} If $h\in\cQ$ and $\sigma(h)\subset I$, then $f(h)\in\cS$.
\itemitem{(ii)} If $p, t\in B(H)_{\text sa}$ such that $p$ is a projection and $\sigma(t)\subset I$, then $pf(ptp)p\le f(t)$.
\itemitem{(iii)} If $p, t\in B(H)_{\text sa}$ such that  $0\le p\le\b1$ and $\sigma(t)\subset I$ then $f(ptp)\le f(t)+ f(0)(\b1-p)$. 	
\itemitem{(iv)} Either $f=0$ or $f(x)>0$, $\forall x\in I$, and $-1/f$ is operator convex.
\itemitem{(v)} $f$ has a representation.

$$f(x)=c+\int_{r<I}{1\over x-r} d_{\mu_-}(r)+\int_{r>I}{1\over r-x} d_{\mu_+}(r),\leqno(3)$$
where $\mu_{\pm}$ are positive measures such that $\int{1\over 1+|r|} d_{\mu_\pm}(r)<\infty$ and $c\ge 0$.

\demo{proof}\ (i) $\Rightarrow$ (ii): As in the proof of Theorem 3.1 we first prove a matrix version, and we use the same choices of $h$ and the same notation as in the second part of the proof of 3.1. It is no longer true that $f(0)=0$, but since $f(x)\b1=f(x\b1)\in\cS$ for $x\in I$, (i) implies that $f\ge 0$ on $I$. Let $f(a)=a^{''}$ and $k=\sum a^{''}_{ij}e_i \times e_j$. Then $k$ is a compact operator and $k\le s_\infty$, whence $\forall \ep>0$, $k\le s_n+\ep\b1$ for $n$ sufficiently large. For any $n$, the last relation amounts to the matrix inequality $\pmatrix f(a) & 0\\ 0 & 0\endpmatrix \le\pmatrix a'+\ep\b1_k &b'\\ b^{'*} & c'+\ep\b1_l\endpmatrix$. Since $\ep$ is arbitrary, we conclude that $\pmatrix f(a) & 0\\ 0 & 0\endpmatrix \le f(t)$, the matrix version of (ii), the general version follows from the matrix version just as in the proof of 3.1.

\noindent
(ii) $\Rightarrow$ (i): By applying (ii) with $t=x\b1, x\in I$, we deduce $f\ge 0$.  Now let $(p_m)$ be an increasing sequence of finite rank projections in $B(\ell^2)$ which converges weakly (and strongly) to $\b1$ and let $h\in\cQ$ with $\sigma(h)\subset I$. Apply (ii) with $p_m$ for $p$ and $p_{m+1} h_\infty p_{m+1}$ for $t$ to deduce $p_m f(p_m h_\infty p_m)p_m\le f(p_{m+1}h_\infty p_{m+1})$. Multiplying on both sides with $p_{m+1}$, we find $k_m=p_mf(p_mh_\infty p_m)p_m\le p_{m+1}f(p_{m+1}h_\infty p_{m+1})p_{m+1}=k_{m+1}$. Since $k_m\to f(h_\infty)$, it is sufficient to show that $\forall m, \forall\ep>0, k_m\le f(h_n)+\ep\b1$ for $n$ sufficiently large. (Condition (i) of Definition 2.3 follows from the fact that $f\ge 0$.) Since $h_n\to h_\infty$ weakly, $p_mh_np_m\to p_mh_\infty p_m$ in norm, and hence $p_m f (p_mh_np_m)p_m\to k_m$ in norm. Thus for sufficiently large $n$, $k_m\le p_mf(p_mh_np_m)p_m+\ep\b1\le f(h_n)+\ep\b1$. 

\noindent
(i) and (ii) $\Rightarrow$ (iii): Let $0\le p\le\b1$ and choose $k$ in $\cK_{\text sa}$ with $\sigma(k)\subset I$. Choose a sequence $(p_n)$ of projections such that $p_n\to p$ weakly.  (The possibility of this was proved by P. Halmos in[Hal].) Then define $h=\{h_n\}$ in $\cQ$ by $h_n=p_nkp_n, n<\infty$, and $h_\infty=pkp$. Since $\sigma(h)\subset I$, (i) implies $f(h)\in\cS\subset \cW$. So if $t$ is a weak cluster point of $(f(h_n))$, then $t\ge f(pkp)$. But by (ii) and the fact that $f(p_nkp_n)=p_nf(p_nkp_n)p_n+f(0)(\b1-p_n)$, $f(h_n)\le f(k)+f(0)(\b1-p_n)$, and  $f(k)+f(0)(\b1-p_n)\to f(k)+f(0)(\b1-p)$ weakly. Thus $f(pkp)\le t\le f(k)+f(0)(\b1-p)$. As above, this is sufficient to establish the general case of (iii).

\noindent
(iii) $\Rightarrow$ (ii): Apply (iii) with $p$ a projection, and note that $f(ptp)-f(0)(\b1-p)=pf(ptp)p$. 

\noindent
(i) $\Rightarrow$ (iv): We have already seen that (i) implies $f\ge 0$. If $f\ne 0$, let $J$ be an open subinterval of $I$ such that $f(x)>0, \forall x\in J$. If $h\in \cQ$ and $\sigma(h)\subset J$, then $f(h)\in \cS$ and Proposition 2.8 implies $-f(h)^{-1}\in\cW$. Thus Theorem 3.1 implies that $-1/f$ is operator convex on $J$. In particular $-1/f$ is convex, and a convex function cannot approach $-\infty$ at a finite endpoint of its interval of definition. Therefore if either endpoint of $J$ is in $I$, then $f$ does not vanish at that endpoint.

Now let $J_0=\{x\in I^0:f(x)>0\}$, where $I^0$ is the interior of $I$. Then $J_0$ is the disjoint union of open intervals, and the above implies that none of these intervals can have an endpoint in $I^0$. It follows that $J_0=I^0$, and another application of the above shows that $f(x)>0, \forall x\in I$. Now it is clear that $-1/f$ is operator convex on all of $I$.

\noindent
(iv) $\Rightarrow$ (v): We may assume $0\in I^0$, since neither (iv) nor (v) is affected by a translation of the independent variable. Assume $f\ne 0$. Let $\varphi(x)=-1/x$ for $x<0$. Since $\varphi$ is both operator monotone and operator convex, and since $-1/f$ is operator convex, then $f=\varphi(-1/f)$ is operator convex. Let $f$ be represented as in (4) with $x_0=0$. The analytic extension of $f$ into the non--real part of the complex plane is obtained by replacing $x$ by $z$ in (4).

Our first task is to show that $\int 1/|r| d\mu_\pm(r)<\infty$. Let $g(x)=(f(x)-f(0))/x$. One of the main results of J. Bendat and S. Sherman, [BS, Theorem 3.2], implies that $g$ is operator monotone, and also ${g(x)\over f(0)f(x)}={-{1\over f(x)}-(-{1\over f(0)})\over x}$ is operator monotone. By L\"owner's theorem each of these operator monotone functions is either a constant or it carries the upper half plane into itself. If $g$ were a non-zero, constant, then $-1/f$ would not be operator convex; and if $g=0$, then $f$ is a positive constant, a trivial case of (3). And if $g(x)/f(0)f(x)$ is a constant, then $f$ is the reciprocal of a linear function, another trivial case of (3).

Thus we may assume both $g$ and $g/f$ carry the upper half plane into itself (since $f(0)>0$). If $\im z=y>0$, then 
$$\im \ {g(z)\over f(z)}= \im\ {g(z)\over f(0)+zg(z)}={f(0)\im g(z)-y|g(z)|^2\over \text{positive}}.$$
So $f(0)\im g(z)>y|g(z)|^2\ge y|\im g(z)|^2$, whence $\im g(z)< f(0)/y$. 

\noindent
From (4) we obtain
$$g(x)=\int_{r>I}\left({1\over r-x}-{1\over r}\right){1\over r}d\mu_+(r)-\int_{r<I}\left({1\over x-r}-{1\over |r|}\right){1\over |r|}d\mu_-(r)+ax+b,\quad \text{and} \leqno(5)$$

$$\im g(z)=\int_{r>I}{y\over|r-z|^2}{1\over r}d\mu_+(r)+\int_{r<I}{y\over|z-r|^2}{1\over |r|}d\mu_-(r)+ay. \leqno(6)$$

\noindent
This implies $ay^2+\int_{r>I}{y^2\over |r-z|^2}{1\over r}d\mu_+(r)+\int_{r<I}{y^2\over |z-r|^2}{1\over |r|}d\mu_-(r)<f(0)$. If $\re z=0$, then $|r-z|^2=r^2+y^2$, and we can let $y\to \infty$ and apply the monotone convergence theorem to the two integrals. We conclude that $a=0$ and $\int{1\over |r|} d\mu_{\pm}(r)<\infty$. 

Now as mentioned in the proof of Theorem 3.1, the integrands in (4) are $\pm({1\over r-x}-{1\over r}-{x\over r^2})$. We now know that each of these three terms is integrble, so we can drop the linear terms from the integrals and absorb the integrals of the linear terms into $bx+c$, obtaining
$$f(x)=\int_{r>I}{1\over r-x} d\mu_+(r)+\int_{r<I}{1\over x-r} d\mu_-(r)+bx+c, \leqno(4')$$
where $b$ and $c$ no longer have the same values as in (4),
$$g(x)=\int_{r>I}{1\over r(r-x)}d\mu_+(r)-\int_{r<I}{1\over |r|(x-r)}d\mu_-(r)+b\leqno(5')$$
$$\im g(z)=\int_{r>I}{y\over r|r-z|^2}d\mu_+(r)+\int_{r<I}{y\over |r||z-r|^2} d\mu_-(r).\leqno(6')$$
Then the inequality $f(0) \im g(z)>y|g(z)|^2$ yields 
$$\align
&f(0)\bigg(\int_{r>I}{1\over r|r-z|^2}d\mu_+(r)+\int_{r<I}{1\over |r||z-r|^2} d\mu_-(r)\bigg)\\
&> \bigg|b+\int_{r>I}{1\over r(r-z)}d\mu_+(r)-\int_{r<I}{1\over |r|(r+z)}d\mu_-(r)\bigg|^2.\endalign$$
If we let $z\to \infty$ so that $\im z$ is bounded away from 0, the dominated convergence theorem applies and gives $f(0)\cdot 0\ge |b|^2$, whence $b=0$. 

Now we calculate $\lim\limits \Sb y\to\infty\\ \re z=0\endSb y\im g(z)=\lim\limits_{y\to\infty} {y^2\over r^2+y^2}{1\over |r|} d\mu(r)$, where $\mu=\mu_++\mu_-$. The monotone convergence theorem applies and yields $\lim y \im g(z)=\int{1\over |r|}d\mu(r)$. Since $\im g(z)<f(0)/y$, we conclude $\int{1\over |r|}d\mu(r)\le f(0)=c+\int{1\over |r|} d\mu(r)$, whence $c\ge 0$.

\noindent
(v) $\Rightarrow$ (i): As in the proof of Theorem 3.1, we can calculate $f(h)$, for $h\in\cQ$ and $\sigma(h)\subset I$, by substituting $h$ for $x$ in (3), thus obtaining a Bochner integral. Because of the properties of $\cS$ established in Proposition 2.7, it is enough to show that each value of the integrand is in $\cS$ and observe that $c\b1\in\cS$.  The first fact follows from Proposition 2.8.
\enddemo

\example{Remarks 3.3} (i) The assumption that $0\in I$ was made because condition (iii) doesn't make sense otherwise. But conditions (i), (iv), and (v) all make sense for arbitrary intervals $I$ and are unaffected by translations of the independent variable. Also, condition (ii) can easily be interpreted to make sense for arbitrary $I$, and then it too is unaffected by translation. For example, we can extend $f$ arbitrarily to $I\cup \{0\}$, define $f(ptp)$ by the Borel functional calculus, and observe that $pf(ptp)p$ depends only on $f_{|I}$.

Thus for $f$ defined on an arbitrary interval $I$, we can define strong operator convexity by applying any of conditions (i), (ii), (iv), or (v) directly to $f$, or by applying condition (iii) to $f(\cdot +x_0)$ for some $x_0$ in $I$. The last yields

\noindent
(iii)$_{x_0}$ $f(ptp+x_0(\b1-p^2))\le f(t)+f(x_0)(\b1-p)$, for $p, t \in B(H)_{\text sa}$, $0\le p\le\b1$, and $\sigma(t)\subset I$. 

So what we have proved shows that (iii)$_{x_0}$ is independent of the choice of $x_0$ in $I$.

(ii) If $f$ is real analytic on $I^0$ and strongly operator convex on some non-empty open subinterval $J$ of $I$, then $f$ (still assumed continuous and real--valued on $I$) is strongly operator convex on $I$. This can be proved easily from condition (iv) (also easily from (v)). The assumptions and the principle of uniqueness of analytic continuation imply that there is a holomorphic function $\tf$ on $I^0\cup\{z\in\bC: \im z\ne 0\}$ which agrees with $f$ on $I^0$. So (iv) for $f, J$ implies that $\tf(z)\ne 0$ for $\im z>0$ and $ {-{1\over\tf(z)}-\left(-{1\over\tf(x_0)}\right)\over z-x_0}$ is either constant or maps the upper half plane into itself. And this fact about $\tf$ implies that $f$ is strongly operator convex on any open subinterval $J$, of $I$ which contains $x_0$ and does not contain any zeros of $f$. Then part of the proof above that (i) $\Rightarrow$ (iv) applies to show $f(x)>0$, $\forall x\in I$, whence the conclusion.

(iii) If an operator monotone or operator convex function $f$ on an open interval $I$ is extended to one or both endpoints of $I$ so that it is still monotone or convex but no longer continuous, then the operator inequalities used to define operator monotonicity or operator convexity will still hold for the extended function (using the Borel functional calculus to define $f(t_i)$). Let $f$ be the function on $[0,\infty)$ with $f(0)=1$ and $f(x)=0$ for $x>0$. Then for $t\ge 0$ in $B(H)$, $f(t)$ is the kernel projection of $t$. Then $f$ satisties conditions (ii) and (iii) of the theorem but fails conditions (i), (iv), and (v). With regard to (i), if $h$ is a positive element of $\cQ$, $f(h)$ need not even be in $\cW$. The easist way to prove (ii) and (iii) for $f$ is to note that $f(x)=\lim\limits_{\ep\to 0^+}\ep/(\ep+x)$ and the function $f_\ep(x)=\ep/(\ep+x)$ is strongly operator convex on $[0,\infty)$. The dominated convergence theorem and the spectral theorem imply that $f_\ep(t)\to f(t)$ strongly. For the function $1+f$, (ii), (iii), and (iv) hold and (i) and (v) still fail.

(iv) It is easy to construct strongly operator convex functions from operator convex functions with the help of condition (iv). If $g$ is operator convex on an interval $I$, choose $\lambda\in\bR$ such that $g(x)<\lambda$ for some $x$ in $I$. Then let $J$ be a subinterval of $I$ such that $g(x)<\lambda$, $\forall x\in J$, and let $f=1/(\lambda-g_{|J})$. 
\endexample

\subhead 4. \ Additional operator inequalities\endsubhead

Part of the original motivation for condition (iii) in Theorem 3.2 (and [B1, Theorem 2.36]) was to obtain a condition which is related to (ii) in the same way as (2) relates to (1). However, it is not at all clear that (iii) accomplishes this. Although it is obvious that (iii)$\Rightarrow$ (ii)$\Rightarrow$ (1), it is not obvious that (iii)$\Rightarrow$ (2). The only way we know to prove this is to repeat part of the proof of [HP, Theorem 2.1] and deduce (2) from (1). We do not understand how (iii) fits into the general scheme of things and note that it is a somewhat peculiar looking  condition. (We did make direct use of (iii) in [B1] in the proof that (iii)$\Rightarrow$ (iv) in Theorem 2.36, but this was not a true application of (iii), since we could have proved (iv) just as efficiently in a different way.) And condition (iii)$_{x_0}$ in Remark 3.3 (i) looks even more peculiar.

In this section we take a different approach to obtain some operator inequalities for strongly operator convex functions which are similar to but stronger than some operator inequalities for (general) operator convex functions. The idea is to find a way to deduce an inequality from (1), and then see what stronger inequality we get if we use 3.2 (ii) instead of (1). We always exclude the case $f=0$ in what follows.

The following inequality, due to Hansen and Pedersen  though not quite explicity stated in [HP], holds if $f$ is operator convex on $I$, $a_1, \dots, a_n\in B(H)$, $t_1, \dots, t_n\in B(H)_{\text {sa}}$, $\sigma(t_i)\subset I$, and $\sum a_i^*a_i=\b1$. 

\noindent
(7)\quad $f(\sum a_i^* t_i a_i)\le\sum a_i^*f(t_i)a_i$.

\noindent
To deduce (7) from (1) consider the isometry from $H$ into $H\oplus\cdots\oplus H$ given by the column $v=\pmatrix a_1\\ \vdots\\ a_n \endpmatrix$. If the range of $v$ is $M$, let $H'$ be a Hilbert space of the same dimension as $M^\bot$  and find an isometry $w=\pmatrix b_1\\ \vdots\\ b_n \endpmatrix$ from $H'$ onto $M^\bot$. Thus $u=(v\  w)$ is a unitary from $H\oplus H'$ to $H\oplus\cdots\oplus H$. Then (7) results from applying (1) to $t=u^*(t_1\oplus\dots\oplus t_n)u\in B(H\oplus H')_{\text {sa}}$ and $p$ the projection with range $H$. If $f$ is strongly operator convex and we instead apply 3.2 (ii) with the same data, we obtain the following, where the $2\times 2$ matrices represent elements of $B(H\oplus H')$. 

\noindent
(8)\quad $\pmatrix f(\sum a_i^* t_i a_i) & 0\\ \   \\ 0 & 0 \endpmatrix\le \pmatrix \sum a_i^* f(t_i) a_i & \sum a_i^* f(t_i) b_i\\ \   \\  \sum b_i^* f(t_i) a_i & \sum b_i^* f(t_i) b_i \endpmatrix$. 

\noindent
Obviously (8) implies (7), and we can deduce an inequality in $B(H)$ by applying the following principle. 

\noindent
(9) $\pmatrix a &b\\b^* &c\endpmatrix\ge 0 \Leftrightarrow a\ge bc^{-1}b^*$, provided $c$ is positive and invertible.

\noindent
Since $f(x)>0, \forall x\in I$, (9) does apply in our situation, and we obtain

\noindent
(10) $f(\sum a_i^*t_ia_i)\le \sum a_i^*f(t_i)a_i-(\sum a_i^*f(t_i)b_i)(\sum b_i^*f(t_i)b_i)^{-1}(\sum b_i^*f(t_i)a_i)$.

\noindent
Of course (8) and (10) are not fully explicit because we haven't given formulas for $b_1,\dots, b_n$, but this could be remedied, at the cost of more complicated notation, as follows: The projection $q$ with range $M^\bot$ is given by the $n\times n$ operator matrix $\b1_n-vv^*$. Let $H'=M^\bot$ and let $w$ be the inclusion of $M^\bot$ into $H\oplus\cdots\oplus H$. It would actually be easiest, and permissible, to replace $w$ by $q$, so that $u$ becomes  a co-isometry given by an $n\times(n+1)$ matrix with entries in $B(H)$. 

\noindent
We prefer instead to stick with (8) and (10) and will consider a couple of special cases where $b_i, \dots, b_n$ can be calculated more easily. Since strongly operator convex functions do not satisfy the condition $f(0)\le 0$, replace (2) by:

\noindent
(2$'$)\  $f(a^*ta)\le a^*f(t)a+f(0)(\b1-a^*a)$. 

\noindent
This is the special case of (7) where $n=2, t_2=0$, and $a_2=(\b1-a^*_1a_1)^{\frac 12}=(\b1-a^*a)^{\frac 12}$. Then, by a well known formula, we can take $H'=H$, $b_1=(\b1-aa^*)^{\frac 12}$, and $b_2=-a^*$. Thus when $f$ is strongly operator convex, we get the following strengthening of (2$'$);

\noindent
(2$''$)\ $f(a^*ta)\le a^*f(t)a+f(0)(\b1-a^*a)-a^*(f(t)-f(0)\b1)(\b1-aa^*)^{\frac 12}((\b1-aa^*)^{\frac 12}f(t)(\b1-aa^*)^{\frac 12}+f(0)aa^*)^{-1}(\b1-aa^*)^{\frac 12}(f(t)-f(0)\b1)a$.

\noindent
A better comparison with (2) can be obtained by rewriting (2$''$) in terms of $g=f-f(0)$. (Here $f$ is still strongly operator convex but $g$ isn't.)

\noindent
(2$'''$) \ $g(a^*ta)\le a^*g(t)a-a^*g(t)(\b1-aa^*)^{\frac 12}(f(0)\b1+(\b1-aa^*)^{\frac 12} g(t)(\b1-aa^*)^{\frac 12})^{-1}(\b1-aa^*)^{\frac 12}g(t)a$.

\noindent
If $a$ is a projection, then (2$''$) implies 3.2 (ii), whence (2$''$) is equivalent to strong operator convexity.

Next take $n=2$, $a_1=\lambda^{\frac 12}\b1$, and $a_2=(\b1-\lambda)^{\frac 12}\b1$ for $\lambda\in (0,1)$, so that (7) becomes the defining relation for operator convexity. Since $a_2\ge 0$, as above there are simple formulas, $b_1=(1-\lambda)^{\frac 12}\b1$ and $b_2= -\lambda^{\frac 12}\b1$. Then for $f$ strongly operator convex (10) becomes:

\noindent
(11)\ $$\align
f(\lambda t_1+(1-\lambda)t_2)&\le \lambda f(t_1)+(1-\lambda)f(t_2)\\
&-\lambda (1-\lambda)(f(t_1)-f(t_2))((1-\lambda)f(t_1)+\lambda f(t_2))^{-1}(f(t_1)-f(t_2)).\endalign$$

\noindent
Note that (11) is not satisfied by arbitrary operator convex functions, for example by $f(x)=x^2$, even if $t_1$ and $t_2$ are scalar operators.

\noindent
By using the same construction that Davis used in [D1] to show that operator convexity  implies (1), we can show that (11) implies 3.2 (ii). So (11) is equivalent to strong operator convexity. 

Next we consider the case where the subspace $M^\bot$ of $H\oplus\dots\oplus H$ is one--dimensional. This implies, if $H=\ell^2$, that the $C^*$--algebra generated by $a_1\dots, a_n$ is an extension of $\cK$ by the Cuntz algebra $\Cal O_n$. This extension does not represent the usual generator of Ext$_s(\Cal O_n)$, the group of extensions of $\cK$ by $\Cal O_n$ with strong equivalence, but rather its negative. Let $v=u_1\oplus\cdots\oplus u_n$ be a unit vector in $M^\bot$. Then for $f$ strongly operator convex (10) becomes:

\noindent
(12)\ $$\align
f\left(\sum a_i^*t_ia_i\right)&\le \sum a_i^*f(t_i)a_i\\
&-\left(\sum(f(t_i)u_i, u_i)\right)^{-1}\left(\sum a^*_if(t_i)u_i\right)\times \left(\sum a_i^*f(t_i)u_i\right).\endalign$$

\noindent
Finally, consider the special case of 3.2 (ii) where the projection $p$ has rank one. (This is not an additional inequality in the sense meant by the title of this section.) This amounts to the following situation: We are given a probability measure $\mu$ supported on a compact subset of $I$, $H=L^2(\mu)$,  $t$ is multiplication by the identity function on $I$, and the range of $p$ is the set of constant functions in $L^2$. In this case (1) becomes:

\noindent
(13)\ $f(\int x d\mu(x))\le \int f(x)d\mu(x)$. 

\noindent
Of course (13) is just the classical Jensen's inequality, which is valid for arbitrary convex functions. But 3.2 (ii) yields, when $f$ is strongly operator convex:

\noindent
(14) \ $f(\int x d\mu(x))\le 1/\int f(x)^{-1}d\mu(x)$.

\noindent
Also 3.3(iii) yields

\noindent
(15) \ $f(s^2 \int d\mu(x))\le 1/\int f(x)^{-1}d\mu(x) +f(0)(1-s)$, $0\le s\le 1$.

\noindent
That (14) implies (13) is just the fact that the harmonic mean is less than or equal to the arithmetic mean. And (14) is not true for arbitrary operator convex functions. We have not found any applications of (14), only new proofs of already known facts, but conceivably (14) could have an interesting consequence if applied to a particularly interesting strongly operator convex function $f$. If the measure $\mu$ is supported by a two--point set, then (14) is the same as the specialization of (11) to scalar operators.  And we have found no interesting consequences of (15).

\subhead 5. \ A differential criterion\endsubhead

If $f$ is a smooth function on an open interval $I$, then the function $t\mapsto f(t)$ for self--adjoint $n\times n$ matrices $t$ with $\sigma(t)\subset I$ is also smooth. In this section we will denote this function on $\m$ by $F_n$. (Up to now we have been casual about the notation.) The first derivative of $F_n$ at $t$ is a linear function from $\m$ to itself, $h\mapsto F'_n(t)\cdot h$, and the second derivative is a symmetric bilinear function from $\m\times \m$ to $\m$, $(h, k)\mapsto F''_n(t)(h, k)$. A well known criterion for operator monotonicity is that $F'_n(t)\cdot h\ge 0$ whenever $h\ge 0$, for arbitrary $n$. And a well known criterion for operator convexity is that $F''_n(t)(h, h)\ge 0$, for arbitrary $n$. Of course, these criteria aren't complete without information on how to compute $F'_n$ and $F''_n$ in terms of $f$, and the reader is referred to the existing literature for this.  Condition (iv) of Theorem 3.2 true makes it easy to derive a differential criterion for strong operator convexity. It is the following $2n\times 2n$ matrix inequality, which has to hold for arbitrary $n$:

\noindent
(15)\ $$\pmatrix F^{''}_n(t)(h, h)/2 & F'_n(t)\cdot h\\
\     \\
 F'_n(t)\cdot h & F_n(t)\endpmatrix\ge 0.$$

\proclaim{Theorem 5.1}\ If $f$ is a continuous real--valued function on an interval $I$ which is $C^2$ on $I^0$, then $f$ is strongly operator convex if and only if (15) holds for all $n$, for all $t$ in $\m$ with $\sigma(t)\subset I^0$, and for all $h$ in $\m$. 
\endproclaim

\demo{Proof}\ If $J$ is any open subinterval of $I$ such that $f(x)>0$, $\forall x\in J$, then $f_{|J}$ is strongly operator convex if and only if $
G_n^{''}(t)(h,h)\ge 0$  for arbitrary $n$. Here $g=-1/f_{|J}$ and $G_n$ relates to $g$ as $F_n$ to $f$. Computation shows that 
$G'_n(t)\cdot h=F_n(t)^{-1}(F'_n(t)\cdot h) F_n(t)^{-1}$, and 
$$\align
&G^{''}_n(t)(h, k)=-F_n(t)^{-1}(F'_n(t)\cdot k)F_n(t)^{-1}(F'_n(t)\cdot h)F_n(t)^{-1}\\
&+F_n(t)^{-1}(F^{''}_n(t)(h, k))F_n(t)^{-1}-
F_n(t)^{-1}(F'_n(t)\cdot h) F_n(t)^{-1}(F'_n(t)\cdot k)F_n(t)^{-1}.\endalign$$
Thus
 $$\align
 G^{''}_n(t)(h, h)&=F_n(t)^{-1}(F^{''}_n(t)(h, h))F_n(t)^{-1}\\
 &-2 F_n(t)^{-1}(F'_n(t)\cdot h) F_n(t)^{-1}(F'_n(t)\cdot h) F'_n(t)^{-1}.\endalign$$
  So $G^{''}_n(t)(h,h)\ge 0$ if and only if $F^{''}_n(t)(h, h)\ge 2(F'_n(t)\cdot h) F_n(t)^{-1}(F'_n(t)\cdot h)$. This is equivalent to (15), for $\sigma(t)\subset J$, by (9). 

Now it is clear that if $f$ is strongly operator convex on $I$, which implies that it is real analytic on $I^0$, then (15) holds. Conversely, if (15) holds, then part of the proof that (i) $\Rightarrow$ (iv) in Theorem 3.2 shows that $f(x)>0$, $\forall x\in I$, if $f$ is not identically $0$.  (Note that the lower righthand corner of (15) implies $f\ge 0$.) Thus $f$ is strongly operator convex on $I^0$ and also on $I$.
\enddemo

\example{Remark}  The existing literature on operator convexity shows that it is not necessary to prove $F_n''(t)(h,h) \ge 0$ for all pairs $(t,h)$ but only for certain well chosen pairs.  The same applies to Theorem 5.1.
\endexample

\bye